\documentclass{article}

\usepackage{amssymb}

\usepackage{epsfig}

\newcommand{\actionfunctional}{{\cal{A}} _{H, \sigma}}

\newcommand{\comphalfaab}{CF ^{[a,b)}(H;\alpha)}

\newcommand{\cotgb}{T^{\ast}M}

\newcommand{\eseuno}{S^{1}}

\newcommand{\fhabalpha}{HF ^{[a,b)}(H;\alpha)}

\newcommand{\kplus}{K^{+}}

\newcommand{\kminus}{K^{-}}

\newcommand{\lsub}{l_{\gamma _{\alpha}}} 

\newcommand{\lsubi}{l_{\gamma _{\alpha _i}}}

\newcommand{\perhalfaab}{{\cal{P}} ^{[a,b)}(H; \alpha)}

\newcommand{\piunotilde}{\widetilde{\pi_{1}}(M)}

\newcommand{\reals}{\mathbb{R}}

\newcommand{\integers}{\mathbb{Z}}

\newcommand{\floermag}{HF^{[a,b)} _{\omega}}

\newcommand{\floerst}{HF^{[a,b)}_{\omega_{0}}}

\newtheorem{prop}{\bf{Proposition}}

\newtheorem{lemma}{\bf{Lemma}}

\newtheorem{teo}{\bf{Theorem}}

\newtheorem{cor}{\bf{Corollary}}

\title{Non-contractible periodic orbits of Hamiltonian flows on twisted cotangent bundles}

\author{C\'esar J. Niche \thanks{This work was partially supported by the NSF. The author acknowledges support from IMERL, Fac. de Ingenier\'{\i}a - UDELAR, Montevideo, 
Uruguay.}\\ Mathematics Department \\ University of California at Santa Cruz \\ Santa Cruz CA 95064 \\ USA \\ E-mail: cniche@math.ucsc.edu}

\begin{document}

\noindent

\maketitle

\begin{abstract}
For many classes of symplectic manifolds, the Hamiltonian flow of a function with sufficiently large variation must have a fast periodic orbit. This principle is the base of the notion of Hofer-Zehnder capacity and some other symplectic invariants and leads to numerous results concerning existence of periodic orbits of Hamiltonian flows. Along these lines, we show that given a negatively curved manifold $M$, a neigbourhood $U_{R}$ of $M$ in $\cotgb$, a sufficiently $C^{1}$-small magnetic field $\sigma$ and a non-trivial free homotopy class of loops $\alpha$, then the magnetic flow of certain Hamiltonians supported in $U_{R}$ with big enough minimum, has a one-periodic orbit in $\alpha$. As a consequence, we obtain estimates for the relative Hofer-Zehnder capacity and the Biran-Polterovich-Salamon capacity of a neighbourhood of $M$.

\end{abstract}

\section{Introduction}

We address the existence problem for periodic orbits of Hamiltonian flows on twisted cotangent bundles. These flows model the dynamics of a charged particle in a magnetic field. More specifically, we focus here on the case where the base is equipped with a metric of negative scalar curvature. Let us start our discussion by describing some recent results which are relevant to us. \\

It is a well-known fact that the flow generated by a $C^{2}$-small Hamiltonian has no non-trivial one-periodic orbits, i.e., all its periodic orbits are fixed points. It follows then, that for the flow to have such periodic orbits, the Hamiltonian must have, in some sense, a big enough oscillation. The next result is a particular instance of this general principle.

\begin{teo} {\bf (Thm 2.4, \cite{gg})} 
\label{teo-1} Let $M$ be a closed symplectic submanifold of a geometrically bounded, symplectically aspherical manifold $W$ and let $U$ be a sufficiently small neighbourhood  of $M$. Let ${\cal{H}}(U,M)$ be the collection of compactly supported Hamiltonians that are constant in a neighbourhood of $M$ and such that $\max _{U} H = H \vert _{M}$. Then there exists a constant $C = C(U) > 0$ such that the Hamiltonian flow of every $H \in {\cal{H}}(U,M)$ with $\max H > C$ has a non-trivial contractible one-periodic orbit with action greater than $\max H$.
\end{teo}

Some of the key ideas in the proof of Theorem \ref{teo-1} are also building blocks of the proof of the following theorem, which was proved for negatively curved manifolds in \cite{bps} and for arbitrary manifolds in \cite{weber}.

\begin{teo}  \label{teo-2}Let $M$ be a closed, connected manifold, $(\cotgb, \omega_0)$ its cotangent bundle  endowed with the canonical symplectic form $\omega_0 = d\lambda$ and $H:\cotgb \to \reals$ a proper, bounded below Hamiltonian. Assume that a sublevel set $\{ H < c \}$ contains the zero section. Then, for every non-trivial  class $\alpha$ in the set of free homotopy classes of loops in $M$, there is a dense set of energy values $e \in (c, + \infty)$ such that the level set $\{H = e\}$ contains a closed orbit $x_e$ in the class $\alpha$, with $\int _{x_{e}} \lambda > 0$. 
\end{teo}

As we said in the previous paragraph, the proofs of these results are related, as they are based on arguments relying on filtered Floer homology, symplectic homology and symplectic capacities. We recall that the filtered Floer homology of a Hamiltonian is the homology of a chain complex generated by one-periodic orbits whose actions lie in a certain range. Non-triviality of this Floer homology would guarantee the existence of one-periodic orbits, but to explicitly compute it from the definition, we need precisely these same periodic orbits to generate our chain complex. This circular argument is a problem already present in Floer's proof of Arnold's conjecture. The basic idea for solving it, is to relate the Floer homology to the geometry and topology of the manifold on which the dynamics takes place. We roughly describe now the way in which this idea is carried out in the articles just mentioned. Given a Hamiltonian $H$, two functions are constructed, bounding $H$ from above and below, respectively. These functions ``capture'' some information about the geometry of the manifold and as a result, their filtered Floer homologies are amenable to computation and are non-trivial. Then, a special monotone homotopy between these squeezing functions is defined and as a consequence of standard results, their filtered Floer homologies are isomorphic. This isomorphism factors through the filtered Floer homology of $H$, which is thus non-trivial. This implies the existence of a one-periodic orbit for the flow of $H$. \\

We introduce now the setting in which we work and state our main results. Let $(M,g)$ be a closed, negatively curved Riemannian manifold and $U_{R}$ a neighbourhood of finite radius $R$ about $M$ in $\cotgb$. Let $\sigma \in \wedge ^{2} (M)$ be a closed $2$-form, which we call the magnetic field, and $\omega_{0} = d \lambda$ the standard symplectic form on $\cotgb$. Let $\alpha$ be a fixed non-trivial class in $\piunotilde$, the space of free homotopy classes of loops in $M$. Throughout the article, we identify $\piunotilde$ and $\widetilde{\pi_{1}}(\cotgb)$.

\begin{teo}
\label{first-theo}
Let $H \in C^{\infty}(U_{R} \times S^{1})$. Then, there exist constants $C = C_{\alpha} (R)$ and $\epsilon = \epsilon (H,g) > 0$ such that if $\min _{M \times S^{1}} H > C$, the flow of $H$ with respect to $\omega = \omega_{0} + \pi ^{\ast} \sigma$ has a one-periodic orbit in $\alpha$, for any magnetic field $\sigma$ with $\Vert \sigma \Vert _{C^{1}} < \epsilon$. 
\end{teo}

We now introduce the class ${\cal{H}} _{V} \subset C^{\infty} (\cotgb \times \eseuno)$ of Hamiltonians $H$ such that

\begin{itemize}
\item [a)] $supp \, H =  U_{R} \times \eseuno$;
\item [b)] $H$ depends only on $t \in S^{1}$ in a fixed neighbourhood $V$ of $M$;
\item [c)] $\max _{U_{R}} H_{t} = H_{t} \vert _{V}$, for all $t \in S^{1}$.
\end{itemize}

\begin{teo}

\label{second-theo}
There exist constants $C = C_{\alpha}(R)$ and $\epsilon = \epsilon (V,g) > 0$ such that for any $H \in {\cal{H}} _{V}$ with $\min _{V \times \eseuno} H > C$, the flow of $H$ with respect to $\omega = \omega_{0} + \pi ^{\ast} \sigma$ has a one-periodic orbit in $\alpha$, for any magnetic field $\sigma$ with $\Vert \sigma \Vert _{C^{1}} < \epsilon$.

\end{teo}

{\bf Remark} The value of the constant $C = C_{\alpha}(R)$ is computed in section \ref{constru-k}, where it is shown that one can take $C = C_{\alpha}(R) = R \, \lsub$, where $\lsub$ is the length of the unique closed geodesic $\gamma$ in the class $\alpha$ for the geodesic flow of $g$. \\

Let $c_{HZ} (\cdot, \cdot)$ be the relative Hofer-Zehner capacity introduced in \cite{gg}. Then, from this remark and Theorem 2.9 from \cite{gg} we obtain the following 

\begin{cor}
\label{first-coro}
\begin{displaymath}
c_{HZ} (U_{R}, V) \leq R \, \inf _{\alpha \in \piunotilde} \lsub.
\end{displaymath}
\end{cor}

The class ${\cal{H}} _{V}$ plays an essential role in Theorem \ref{second-theo}, allowing us to bound the magnitude of the magnetic field {\it in terms of $V$ and $g$ only}. In other words, {\it every} magnetic field with magnitude less than $\epsilon = \epsilon (V,g)$ is such that the magnetic flow of {\it every} Hamiltonian with big enough oscillation has a $1$-periodic orbit in $\alpha$. If the Hamiltonian $H$ is not in ${\cal{H}} _{V}$, Theorem \ref{first-theo} still guarantees the existence of a $1$-periodic orbit in $\alpha$. However, in this case, the magnitude of the magnetic field {\it depends on $H$}.

In the following proposition we construct an example that shows that Theorem \ref{second-theo} fails if $H$ is not in ${\cal{H}} _{V}$. It also shows that the magnetic term in Theorem \ref{first-theo} may be arbitrarily small. We remark that the Hamiltonians constructed in this example can have arbitrarily large oscilations.

\begin{prop} \label{propo-1}
Let $M$ be a surface endowed with a metric of constant negative curvature equal to $-1$ and $C > 0$ an arbitrarily large constant. Then, there exists a sequence of Hamiltonians $\{ K_{n} \} _{n \in \mathbb{N}}$ with compact support in $ U_R \subset \cotgb$, eventually leaving ${\cal{H}} _{V}$ and with $\min _{M \times \eseuno} K_{n} > C$  and a sequence of magnetic fields $\{ \sigma _{n} \} _{n \in \mathbb{N}}$, with $\Vert \sigma _{n} \Vert _{C^{1}} \to 0$, such that the periodic orbits of the Hamiltonian flow of $K_{n}$ with respect to $\omega_{n} = \omega _{0} + \pi ^{\ast} \sigma _{n}$, are all contractible.
\end{prop}

We recall now the definition of the $\pi_{1}$-sensitive relative symplectic capacity introduced in \cite{bps}. Let $(N, \omega)$ be an open symplectic manifold and ${\cal{H}} (N)$ be the family of Hamiltonians with compact support  in $\eseuno \times N$. Let us denote by ${\cal{P}} (H; \alpha)$ the set of $1$-periodic orbits of the flow of $H$ in the class $\alpha$ and by ${\cal{H}} _{k} (N, A)$ the class

\begin{displaymath}
{\cal{H}} _{k} (N, A) = \{ H \in {\cal{H}} (N): \inf _{\eseuno \times A} H \geq k \}
\end{displaymath}

where $A$ is a compact subset of $N$. We then define the relative symplectic capacity

\begin{displaymath}
c_{BPS} (N,A; \alpha) = \inf \{ k > 0: {\cal{P}} (H; \alpha) \neq \emptyset, \forall H \in {\cal{H}} _{k} (N, A) \}.
\end{displaymath}

Similarly to Corollary \ref{first-coro} we obtain

\begin{cor}

\begin{displaymath}
c_{BPS} (U_{R}, V; \alpha) \leq R \lsub.
\end{displaymath}
\end{cor}

When $\sigma = 0$, it follows from Theorem 3.2.1 in \cite{bps} that $c_{BPS} (U_{R}, M;\alpha) = R \lsub$. However, as an immediate consequence of Proposition \ref{propo-1} we obtain 

\begin{cor} \label{propo-2}
For every $\epsilon > 0$, there exists a magnetic field $\sigma$ with $\Vert \sigma  \Vert _{C^{1}} < \epsilon$ such that
\begin{displaymath}
c_{BPS} (U_{R}, M;\alpha) = + \infty.
\end{displaymath}
\end{cor}

This Corollary expresses the failure of Theorem \ref{second-theo} in terms of this special symplectic capacity. \\

There are many results concerning existence of periodic orbits for Hamiltonian flows on twisted cotangent bundles. In the case of non-trivial contractible orbits of twisted geodesic flows, existence has been proved, under various hypothesis in \cite{cgk}, \cite{frau-schlenk}, \cite{gg}, \cite{ginz-kerman-1}, \cite {ginz-kerman-2} and \cite{xuxa-2}. For results about existence of non-contractible orbits, see  \cite{bps},  \cite{gat-lal}, \cite{yi-jen-lee} and \cite{weber}. \\

This article is organized as follows. In section \ref{tools} we introduce all the definitions and statements we need to prove our result. More precisely, in section \ref{neg-curv} we describe the dynamics of geodesic and magnetic flows and state some dynamical properties of the geodesic flows of negatively curved Riemannian manifolds. Then, in section \ref{floer} we define filtered Floer homology and state some of its basic properties. In section \ref{proofs}, we prove Theorem \ref{second-theo}. We omit the proof of Theorem \ref{first-theo} as it  is analogous to the one of Theorem \ref{second-theo}. In section \ref{constru-k} we define the squeezing functions mentioned in the Introduction. We compute their Floer homologies, working with the standard symplectic form in section \ref{compute} and with the twisted symplectic form in section \ref{perturb}, where we complete the proof of the Theorem following the ideas described after Theorems \ref{teo-1} and \ref{teo-2}. Finally, in section \ref{counter}, we prove Proposition \ref{propo-1}.  \\

{\bf Acknowledgments} The author is deeply grateful to Viktor Ginzburg for his constant encouragement, advice and support. He thanks Ba\c sak G\"urel and Ely Kerman for helpful discussions and suggestions.

\section{Tools}
\label{tools}

We give here the definitions and state the basic results we will need. In section \ref{neg-curv}, we introduce the definitions and properties of geodesic and magnetic flows and negatively curved Riemannian manifolds. We follow   \cite{pat} in our presentation of geodesic flows as Hamiltonian ones on $TM$. For the geometry and dynamics of geodesic flows on negatively curved Riemannian manifolds, see   \cite{kling}. Then, in section \ref{floer}, we define the filtered Floer homology, restricted to a non-trivial class $\alpha \in \piunotilde$ of a compactly supported Hamiltonian in a twisted cotangent bundle over a negatively curved Riemannian manifold $M$. As twisted cotangent bundles are geometrically bounded symplectic manifolds (see \cite{audin-laf}, \cite{cgk}), we follow \cite{cgk}. \\

\subsection{Geodesic flows, magnetic flows and negatively curved manifolds}
\label{neg-curv}

Let $(M,g)$ be a closed Riemannian manifold of dimension $n$. The geodesic flow of the metric $g$ is the map $\phi_{t}:TM \to TM$, given by $\phi_{t} (x,v) = (\gamma_{x} (t), \dot{\gamma_{x}}(t))$, where $\gamma_{x}(t)$ is the geodesic such that $\gamma_{x}(0) = x$ and $\dot{\gamma_{x}}(t) = v$. It becomes a Hamiltonian flow for $H = \frac{1}{2} g_{x}(v,v)$, on the symplectic manifold $(TM, g^{*} (\omega _{0}))$, where $\omega _{L} = g^{*} (\omega _{0})$ is the pullback of the standard symplectic form on $\cotgb$ by the metric. 

Using the standard splitting $T_{\theta}TM = H(\theta) \oplus V(\theta)$, where $H(\theta)$ is the lift of $T_{x}M$ by the metric and $V(\theta)$ is the kernel of $(d{\pi}) _{\theta}$, for ${\pi}:TM \to M$, it is easy to check that the geodesic vector field $G:TM \to TTM$  is 

\begin{displaymath}
G(\theta) = (v,0), \quad \theta = (x,v) \in TM.
\end{displaymath}

The isomorphism between $TM$ and $\cotgb$ given by the metric makes the geodesic flow orbit equivalent to the flow of $H = \frac{1}{2} \Vert p \Vert ^{2}$ on $(\cotgb, \omega _{0})$. \\

Let $\sigma \in \wedge ^{2}(M)$ be a closed $2$-form on $M$. The magnetic flow (or twisted geodesic flow) is the Hamiltonian flow of $H = \frac{1}{2} \rho ^2$ on $(TM, \omega _{L} + \pi ^{\ast} \sigma)$. The Lorentz force is the skew-symmetric linear map $Y_{\sigma}: TM \to TM$ such that

\begin{displaymath}
\sigma _{x} (u,v) = g_{x} (Y_{\sigma}(u), v) , \quad u, v \in T_{x}M.
\end{displaymath}

Then the Hamiltonian vector field $X:TM \to TTM $ of $H$ with respect to the form $\omega = \omega_{L} + \pi ^{\ast} \sigma$ can be expressed as 

\begin{displaymath}
X (\theta) = (v, Y_{\sigma}(v))
\end{displaymath}

for $\theta = (x,v) \in TM$. A curve $(\gamma, \dot{\gamma})$ in $TM$ is an integral curve of $X$ if and only if

\begin{displaymath}
\frac{D}{dt} \dot{\gamma} = Y_{\sigma}(\dot{\gamma})
\end{displaymath}

for the covariant derivative $D$. This equation is just Newton's law for a charged particle of unit mass and charge. \\

Negative sectional curvature of the metric $g$ places restrictions on the topology of $M$ and gives the geodesic flow important dynamical properties. We mention here the ones that are relevant to our work. The first property is that every free homotopy class of loops on $M$ contains a unique closed geodesic, which minimizes length in its class. The second one is that, being Anosov, a geodesic flows $\phi _{t}$ on a negatively curved manifold is $C^{1}$ strongly structurally stable. This means that any flow $\psi _{t}: SM \to SM$ which is close enough to $\widetilde{\phi _{t}} = \phi _{t} \vert _{SM}$  in the $C^{1}$ topology, is $C^{0}$ orbit equivalent to it. In other words, there exists a homeomorphism  $h: SM \to SM$, which can be chosen close enough to identity if the perturbation is small enough, such that $ h ^{-1} \, \circ \, \psi _{t} \, \circ \, h$ equals $\widetilde{\phi _{t}}$, up to a time change.

We describe now a property of negatively curved Riemannian manifolds that will allow us to have a well defined action functional for detecting non-contractible one-periodic orbits. We say that a manifold $M$ is {\it atoroidal}, if every continuous map $f: T^{2} \to M$, is such that $f^{\ast}: H^{2}(M; \reals) \to H^{2}(T^{2}; \reals)$ is the zero map. We claim that this is true for a manifold $M$ that admits a metric of negative sectional curvature. We note first that in this case, the induced map $f_{\ast}: \pi _{1} (T^{2}) \to \pi _{1} (M)$ is not a monomorphism, as all Abelian subgroups of $\pi _{1} (M)$ are infinite cyclic (see \cite{kling}, \cite{preis}). After passing to a finite cover $\widetilde{T}^{2}$ and changing coordinates if necessary, we can assume that the lifted map $\tilde{f}:\widetilde{T} ^{2} \to M$ is such that $\tilde{f} _{\ast}(1,0) = 0$, where $\tilde{f} _{\ast}: \pi_{1} (\widetilde{T}^{2}) \to \pi _{1} (M)$ (i.e. the image by $\tilde{f}$ of a meridian $\nu$ of $\widetilde{T}^{2}$, bounds a disk in $M$). The cycle $\tilde{f}(\widetilde{T}^{2})$ is homologous to zero in $M$, as can be seen from the following argument. Let us cut $\widetilde{T}^{2}$ along the curve $\nu$ and close the ends of the cylinder so obtained with disks. This surface is homeomophic to a sphere, so the resulting map $\bar{f}$ can be identified with a map from $S^{2}$ to $M$. As $\pi _{2} (M) = 0$ for a manifold admitting a metric of negative scalar curvature, $\bar{f}$ has to be homotopic to a constant map. As a consequence $\tilde{f}$ extends to the solid torus and hence $[\tilde{f}(\widetilde{T}^{2})]$ is zero in $H_{2}(M;\integers)$. Therefore $[f(\widetilde{T}^{2})]$ is zero in $H_{2}(M;\reals)$. This proves our claim.

\subsection{Filtered Floer homology}
\label{floer}

A twisted cotangent bundle is the symplectic manifold $(\cotgb, \omega = \omega_{0} + \pi ^{\ast} \sigma)$, where $\omega _0 = d \lambda$ is the standard symplectic form on $\cotgb$, and $\sigma \in \wedge ^{2} (M)$ is a closed two form, called the magnetic field. Let ${\cal{L}} \cotgb = C^{\infty}(\eseuno, T^{\ast}M)$ be the loop space of $\cotgb$ and ${\cal{L}} ^{\alpha} \cotgb = \{ x \in {\cal{L}} \cotgb: x \in \alpha \}$, where $\alpha \neq 0$ is a free homotopy class of loops of $\cotgb$. For a neighbourhood $U_{R}$ of radius $R$ about $M$ in $\cotgb$, we consider ${\cal{H}} _{c} = C^{\infty} (U_{R} \times S^{1}, \reals)$, the set of compactly supported, one-periodic Hamiltonians. For a given $H \in {\cal{H}}_{c}$ and a fixed $\gamma \in \alpha$, we define the action functional $\actionfunctional: {\cal{L}} ^{\alpha} \cotgb \to \reals$ by

\begin{equation}
\label{eq:actionfunctional}
\actionfunctional (x) = - \int _{x} \lambda - \int _{A(x,\gamma)} \pi ^{\ast} \sigma + \int _{x} H
\end{equation}

where $A(x,\gamma)$ is the annulus bounded by $x$ and $\gamma$ and $\pi:T^{\ast}M \to M$ is the standard projection. This functional is well defined, as the second term in (\ref{eq:actionfunctional}) does not depend on the annulus $A(x, \gamma)$ due to the atoroidal condition.

It is clear that the critical points of $\actionfunctional$ are one-periodic orbits, in the class $\alpha$, of the Hamiltonian flow of $H$ with respect to the form $\omega = \omega_{0} + \pi ^{\ast} \sigma$. Let us denote  the set of critical points of $(\ref{eq:actionfunctional})$ by  $\cal{P}  (H; \alpha)$ and the subset of such periodic orbits with action between $a$ and $b$ by

\begin{displaymath}
\perhalfaab = \{ x \in {\cal{P}}  (H; \alpha):a \leq \actionfunctional (x) < b \}.
\end{displaymath}

From now on, let us assume that all one-periodic orbits of $H$ are nondegenerate, i.e., $\det (d\phi_{1} ^{H} (x) - Id) \neq 0$, for all $x \in \cal{P}  (H; \alpha)$, where $\phi _{1}$ is the time-one map of the flow of $H$ \footnote{This nondegenracy condition is generic.}.

We then construct the chain complex $(\comphalfaab, \partial)$, where

\begin{displaymath}
\comphalfaab = \bigoplus _{x \in \perhalfaab} \mathbb{Z} _{2} \, x
\end{displaymath}

and $\partial$ is a boundary operator, whose construction we sketch now. 

Let $J_{gb}$ be a fixed almost complex structure for which $(\cotgb, \omega = \omega_{0} + \pi ^{\ast} \sigma)$ is geometrically bounded. Let $\cal{J}$ be the set of smooth $t$-dependent $\omega$-tame almost complex structures that are compatible with $\omega$ near $supp (H)$ and equal to $J_{gb}$ outside a compact set. Each such $J \in \cal{J}$ gives rise to a positive definite bilinear form on the space $T_{x}{\cal{L}} ^{\alpha} \cotgb$. For $x, y \in \perhalfaab$, we consider ${\cal{M}} (x,y,H,J)$, the set of downward gradient trajectories of $\actionfunctional$ with finite energy, that connect $x$ and $y$. For a dense family ${\cal{J}} _{reg} \in {\cal{J}}$, the space${\cal{M}} (x,y,H,J)$ is a finite-dimensional manifold, with 

\begin{displaymath}
dim \, {\cal{M}} (x,y,H,J) = \mu_{CZ} (x) - \mu_{CZ} (y),
\end{displaymath}

where $\mu_{CZ}: \perhalfaab \to \mathbb{N}$ is the Conley-Zehnder index of a non-degenerate periodic orbit. \\

{\bf Remark} We give here a sketch of a proof that the Conley-Zehnder index of $x \in \perhalfaab$ is well defined. Let $\eta$ be a fixed reference curve in $\alpha$, with a fixed reference trivialization of $\eta ^{\ast} T \cotgb$ along $\eta$. Let ${\cal{C}}$ be a cylinder ${\cal{C}}: [0,1] \times \eseuno \to \cotgb$, where ${\cal{C}}(0,t) = \eta (t)$ and ${\cal{C}}(1,t) = x(t)$. We extend the reference trivialization $\eta ^{\ast} T \cotgb$ to ${\cal{C}} ^{\ast} T \cotgb$ and denote the trivialization induced on $x ^{\ast} T \cotgb$ by $\Phi _{x} (t)$. We then use it to define the Conley-Zehnder index in the same way as in the contractible case. We claim that the Conley - Zehnder index of $x$ is independent of $\cal{C}$ (though it depends on the reference trivizalition). Pick a different cylinder ${\cal{C}}'$ and extend $\eta ^{\ast} T \cotgb$ to obtain ${\cal{C}} ^{' \, \ast} T \cotgb$, inducing a new trivialization $\Phi' _{x} (t)$ of $x ^{\ast} T \cotgb$. By gluing ${\cal{C}}$ and ${\cal{C}}'$ appropiately, we obtain a smooth map $f: T^{2} \to \cotgb$, which induces the trivial symplectic pullback bundle $f ^{\ast} (T \cotgb)$ over $T^{2}$, as $c_{1}(f^{\ast}T \cotgb) = 0$ and the dimension of the base is $2$. Without loss of generality we can assume that the trivialization $\Gamma: T ^{2} \times \reals ^{2n} \to f ^{\ast} (T \cotgb)$ is such that it coincides with $\eta ^{\ast} T \cotgb$ along $\eta$. It follows that $\Gamma _{x} (t) = \Gamma \vert _{x(t)}$ is homotopic to $\Phi _{x} (t)$ and $\Phi' _{x} (t)$. This proves that the Conley-Zehnder index of $x$ is independent of the trivialization of $x ^{\ast} T \cotgb$. \\

Assume now that $\mu_{CZ} (x) - \mu_{CZ} (y) = 1$. As $\reals$ acts freely by translation on the time parameter of the gradient flow trajectories, we can define

\begin{displaymath}
\tau (x,y) = \# _{mod 2} \{ {\cal{M}} (x,y,H,J) / \reals \}.
\end{displaymath}

Then, we define the boundary operator 

\begin{displaymath}
\partial (x) = \sum _{y \in \perhalfaab} \tau (x,y) y
\end{displaymath}

where the summation extends over $y$ such that $\mu_{CZ} (x) - \mu_{CZ} (y) = 1$. It can be shown that indeed $\partial ^{2} = 0$.

Then, the Floer homology of $H$ in the interval $[a,b)$ and the class $\alpha$  is the homology of the complex $(\comphalfaab, \partial)$, i.e., 

\begin{displaymath}
\fhabalpha = \frac{Ker \, \partial}{Im \, \partial}.
\end{displaymath}

For details and proofs of these results in the general case of Floer homology, see  \cite{cfh}, \cite{fh} and \cite{salamon}. The proofs for the filtered case are analogous. \\

{\bf Remark} It is unclear whether $\fhabalpha$ depends on $J _{gb}$ or not. It is independent of it, if the set of almost complex structures for which $\cotgb$ is geometrically bounded is connected.  \\ 

We state now some basic lemmas that relate the Floer homology of two functions. Let us denote the action spectrum of $H \in {\cal{H}}$ in $\alpha$ by 

\begin{displaymath}
{\cal S} (H; \alpha) = \{ \actionfunctional (x): x \in {\cal{P}} (H; \alpha) \}.
\end{displaymath}

Let $H$ and $K$ be two functions in 

\begin{displaymath}
{\cal {H}} ^{[a,b)} = \{ H \in C^{\infty} (\cotgb \times \eseuno, \reals) : a, b \notin {\cal S} (H; \alpha) \}
\end{displaymath}

such that $H(x,t) \geq K(x,t)$ in $\cotgb \times \eseuno$. For such functions, there exists a monotone homotopy $s \mapsto K_{s}$ from $H$ to $K$, i.e., a family of functions $K_{s}$ such that 

\begin{displaymath}
K_{s} = \cases{ H \, & $s \in (-\infty, 0]$ \cr 
                K \, & $s \in [1, + \infty)$. \cr } 
\end{displaymath}

and $\partial _{s} K_{s} \leq 0$. This induces a chain complex map $\sigma _{KH} \colon CF ^{[a,b)} (H;\alpha) \to CF ^{[a,b)} (K;\alpha)$ which in turn gives rise to a homomorphism in homology 

\begin{displaymath}
\sigma _{KH}: HF ^{[a,b)} (H;\alpha) \to HF ^{[a,b)} (K;\alpha).
\end{displaymath}

\begin{lemma}
\label{lema-floer-1}
The homomorphism $\sigma _{KH}$ is independent of the choice of monotone homotopy and the following identities hold:

\begin{eqnarray}
\sigma_{KH} \circ \sigma _{HG} & = & \sigma _{KG}, \nonumber \\
\sigma_{HH} & = & id  \nonumber
\end{eqnarray}

for $G \geq H \geq K \in {\cal {H}} ^{[a,b)}$.

\end{lemma}

\begin{lemma}
\label{lema-floer-2}
If $K_{s} \in {\cal {H}} ^{[a,b)}$ for all $s \in [0,1]$, then $\sigma _{KH}$ is an isomorphism.
\end{lemma}

The second lemma says that the only way that $\sigma_{KH}$ can fail to be an isomorphism, is if periodic orbits with action equal to $a$ or $b$ are created during the monotone homotopy. For proofs of these Lemmas, see \cite{cfh}, \cite{fh} and \cite{salamon}. \\

{\bf Remark} We have defined the filtered Floer homology of a Hamiltonian whose $1$-periodioc orbits are nondegenerate. However, it is possible to extend this definition to the space ${\cal {H}} ^{[a,b)}$. Let us endow this set with the strong Whitney $C^{\infty}$ topology. Then, there exists a neighbourhood $V$ containing  $H \in  {\cal {H}} ^{[a,b)}$ such that the filtered Floer homologies of every $K \in  {\cal {H}} ^{[a,b)} \in V$ {\it with nondegenerate periodic orbits} are isomorphic. We can then define $HF ^{[a,b)} (H; \alpha)$, even when $H$ has degenerate periodic orbits, as $HF ^{[a,b)}(K; \alpha)$. For details and proofs, see \cite{bps}, \cite{cfh} and \cite{fh}.

\section{Proof of Theorem \ref{second-theo}}
\label{proofs}

Before proceeding with the proof itself, we give a detailed sketch of it. From now on, we denote by $\floerst$ and $\floermag$ the filtered Floer homologies obtained by using the standard symplectic form $\omega_{0} = d \lambda$ and the twisted symplectic form $\omega = \omega_{0} + \pi ^{\ast} \sigma$ respectively. \\

Let $H \in {\cal{H}} _{V}$. The idea for proving Theorem \ref{second-theo} is to construct functions $\kminus$ and $\kplus$ such that $\kminus \leq H \leq \kplus$ and 

\begin{equation}
\label{eq:isomorphism}
\floermag (\kplus; \alpha) \cong \floermag (\kminus; \alpha) \neq 0
\end{equation}

for some interval $[a,b)$. If we can find a monotone homotopy $\{K_{s}\}$ connecting $\kplus$ and $\kminus$ such that $K_{s} \in {\cal {H}} ^{[a,b)}$, Lemma \ref{lema-floer-2} would imply that $\floermag (H; \alpha) \neq 0$, as the isomorphism (\ref{eq:isomorphism}) factors through this Floer homology group. 

The function $\kplus$ will be constructed in section \ref{constru-k} as an autonomous Hamiltonian whose flow {\it with respect to the standard symplectic form $\omega _{0}$} on every energy level, is a reparametrization of the geodesic flow on that level. To achieve this, we will take $\kplus = \kplus (\rho)$, where $\rho: \cotgb \to \reals$ is the norm on the fibers, i.e. $\rho (x,p) = \Vert p \Vert$. The fact that the geodesic flow of $g$ has a unique closed geodesic in every nontrivial class $\alpha$ when $M$ is negatively curved will force, under certain conditions,  the flow of $\kplus$ with respect to the standard symplectic form to have periodic orbits in $\alpha$. The construction of $\kminus$ will be similar to that of $\kplus$. We remark that both $\kplus$ and $\kminus$ are in  ${\cal {H}} _{V}$. After choosing an appropiate action interval, in section \ref{compute} we will compute $\floerst (\kplus; \alpha)$ and show that it is non-trivial. 

In section \ref{perturb} we prove Theorem \ref{second-theo}. To achieve this, we first prove Proposition \ref{iso-mag-st-kplus}, which says that the Floer homologies of $\kplus$ and $\kminus$ with respect to $\omega$ are non-trivial and isomorphic, provided that the magnetic field $\sigma$ is sufficiently small, where the bound on the magnitude of $\sigma$ can be chosen to depend only on $V$ and $g$. This result follows essentially from the strong structural stability of the geodesic flow and from the fact that $\kplus$ and $\kminus$ are in ${\cal{H}} _{V}$. We then construct a monotone homotopy $\{K_{s} \}$ connecting $\kplus$ and $\kminus$. Proposition \ref{iso-mag-st-kplus} guarantees that the Floer homology of $K_{s}$ remains constant along the homotopy, which proves the theorem.

\subsection{Construction of $K ^{\pm}$}
\label{constru-k}

In this section we construct functions $K^{\pm}$, supported in $U_{R}$ with $\kplus \geq H \geq \kminus$ and such that their flows {\it with respect to the symplectic form $\omega_{0} = d \lambda$} have one-periodic orbits in the non-trivial class $\alpha$. From now on we fix $H \in  {\cal{H}} _{V}$.\\

Before constructing these functions, we find the condition that must hold for their flows to have the required one-periodic orbit in $\alpha$. Let $\kplus = \kplus (\rho)$ be a smooth function supported in $U_{R}$, where $\rho = \rho (x,p) = \Vert p \Vert$ is the norm on the fibers. Note that the flow of $\kplus$ with respect to $\omega_{0}$ is  a reparametrization of the geodesic flow on the energy levels. 

Assume that the loop $x(t) = (q(t), p(t))$ is a critical point of the action functional $A_{\kplus}: {\cal{L}} ^{\alpha} \cotgb \to \reals$ given by 

\begin{equation}
\label{eq:functional-small}
A_{\kplus}(x) = - \int _{x} \lambda + \int _{x} \kplus,
\end{equation}

i.e. $x$ is a one-periodic solution to $\dot{x} = X_{\kplus}(x)$. Then, it is clear that 

\begin{equation}
\label{eq:slope}
(\kplus)' (\rho _{0}) = - \lsub
\end{equation}

where $\Vert p \Vert = \rho _{0}$ and $\lsub$ is the length of the unique closed geodesic in the class $\alpha$. As we want $\kplus \geq H$, we see that the condition $\min _{V \times \eseuno} H > C$, for 

\begin{equation}
C = C_{\alpha} (R) = R \, \lsub
\end{equation}

forces (\ref{eq:slope}) to hold for some $\rho _{0} \in [0, R]$. \\

We now describe $\kplus$ and $\kminus$.  Let $\kplus = \kplus (\rho)$ be a smooth function such that:

\begin{itemize}
\item [a)] $\kplus \geq H$ in $U_{R} \times \eseuno$;
\item [b)] $\kplus$ is even and $\kplus (\rho) = 0$, for $\rho \geq R$;
\item [c)] $\kplus$ is constantly equal to $k^{+} \approx \max _{V \times \eseuno} H$ in an interval $[0, R - \delta)$;
\item [d)] $(\kplus) ' = -\lsub$ only twice in the interval $(R - \delta, R)$. At these points, $(\kplus) '' \neq 0$.
\end{itemize}

\begin{figure}[h]

\begin{center}

\psfig{file=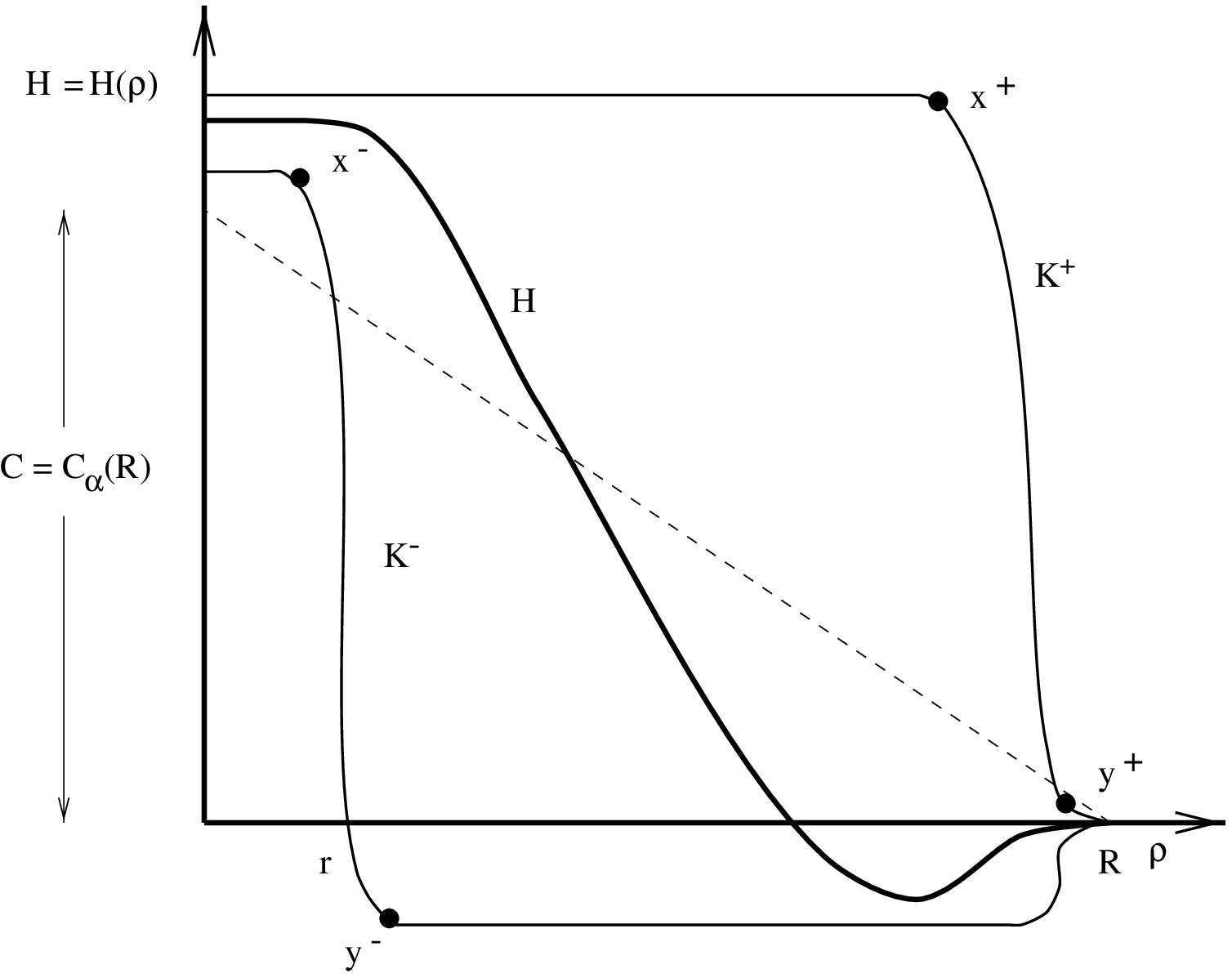, height=9cm, width=6cm, angle=270}

\end{center}

\caption{\label{graficas} Graphs of $\kplus, H, \kminus$. }

\end{figure}

Without loss of generality, we assume that $V$ is a tubular neighbourhood of radius $r_{V}$ about $M$ in $\cotgb$. Let $\kminus$ be the smooth function such that:

\begin{itemize}
\item [a)] $\kminus \leq H$ in $U_{R} \times \eseuno$;
\item [b)] $\kminus$ is even and $\kminus (\rho) = 0$, for $\rho \geq R$;
\item [c)] $\kminus$ is constantly equal to $k^{-}$ in an interval $[0, r)$, with $r > r_{V}$ and $C_{\alpha} (R) < k^{-} < \min _{V \times \eseuno} H$;
\item [d)] $(\kminus) ' = -\lsub$ only twice in an interval $(r, r + \delta ')$. At these points, $(\kminus) '' \neq 0$;
\item [e)] $\kminus$ is constant in $[r + \delta ', R - \delta '')$ and increasing in $[R - \delta '', R]$. 
\end{itemize}

Note that both $\kplus$ and $\kminus$ belong to ${\cal{H}} _{V}$. \\

Let $x^{+}, y^{+}, x^{-}, y^{-}$ be the critical levels where $(K ^{\pm}) ' = -\lsub$ respectively, as seen in Figure \ref{graficas}. The action for the one-periodic orbits on those levels is

\begin{eqnarray}
A_{\kplus}(x^{+}) \approx k^{+} + R \lsub, & \, & A_{\kplus}(y^{+}) \approx R \lsub \nonumber \\ A_{\kminus}(x^{-}) \approx k^{-} + r \lsub, & \, & A_{\kminus}(y^{-}) < A_{\kplus}(y^{+}). \nonumber
\end{eqnarray}

As $k^{-} > C = C_{\alpha}(R)$, by choosing 

\begin{equation}
\label{eq:interval}
R \lsub < a < R \lsub + r \lsub, \qquad b = + \infty
\end{equation}

we can be sure that $x^{+}$ and $x^{-}$ are the only critical levels  with action in $[a,b)$.

\subsection{Computation of Floer homology of $K ^{\pm}$}
\label{compute}

We now compute the filtered Floer homology of $K ^{\pm}$. We closely follow \cite{bps}. \\

We recall that a subset $P \subset \perhalfaab$ is a Morse-Bott submanifold of periodic orbits if $C_{0} = \{x(0): x \in \perhalfaab \}$ is a compact submanifold of $M$ and $T_{x_{0}} C_{0} = \ker (d \phi _{1} (x_{0}) - Id)$ for every $x_{0} \in C_{0}$. The key to the computation of the Floer homology of $K ^{\pm}$ is the following result, which is a version of a theorem from \cite{pozniak} suitable for our setting. \\

\begin{teo} {\bf (\cite{pozniak})} \label{version} Let $H \in {\cal {H}} ^{[a,b)}$. Suppose that the set $P = \{ x \in \perhalfaab: a < A_{H}(x) < b \}$ is a connected Morse-Bott manifold of periodic orbits. Then, $HF ^{[a,b)}(H; \alpha) \cong H _{\ast} (P, {\mathbb{Z}}_{2})$.
\end{teo}

Given $K ^{\pm}$, we denote as ${\cal{P}} _{\omega_{0}}(\rho _{0})$ the subset of ${\cal{P}} ^{[a,b)}(K ^{\pm}; \alpha)$ such that $(K ^{\pm})' (\rho_{0}) = - \lsub$. 

\begin{lemma} {\bf (Lemma 5.3.2, \cite{bps})} \label{lema-uno} The manifold ${\cal{P}} _{\omega_{0}}(\rho _{0})$  is a Morse-Bott manifold of periodic orbits if and only $(K ^{\pm})'' (\rho_{0}) \neq 0$. This manifold is diffeomorphic to $S^{1}$.

\end{lemma}

Then, as a consequence of the construction of $K ^{\pm}$ and of Theorem \ref{version} and Lemma \ref{lema-uno}, we conclude that

\begin{prop} {\bf (\cite{bps})} \label{iso-st-eseuno} Let $[a,b)$ be as in (\ref{eq:interval}). Then 
\begin{displaymath}
HF ^{[a,b)} _{\omega _{0}} (K ^{\pm}; \alpha) \cong H_{\ast} (S^{1}, {\mathbb{Z}}_{2}).
\end{displaymath}

\end{prop}

{\bf Remark} Note that as $K ^{\pm}$ is time independent, the one-periodic orbits of its flow are degenerate, but the Floer homology is still well defined as remarked in section \ref{floer}. \\

\subsection{Proof of Theorem \ref{second-theo}}
\label{perturb}

Let $H \in {\cal{H}} _{V}$ and $K ^{\pm}$ be as in section \ref{constru-k}.

\begin{prop} \label{iso-mag-st-kplus} Let $[a,b)$ be as in (\ref{eq:interval}). Then, there exists $\epsilon = \epsilon (V,g) > 0$ such that 

\begin{displaymath}
\floermag (K ^{\pm}; \alpha) \cong H_{\ast} (S^{1}, {\mathbb{Z}}_{2})
\end{displaymath}

provided $\Vert \sigma \Vert _{C^{1}} < \epsilon$.

\end{prop}

{\bf Proof} From now on, we identify $TM$ and $\cotgb$ through the metric, so as to make use of the structures introduced in section \ref{neg-curv}. We prove the proposition for $\kplus$ only, as the proof for $\kminus$ is analogous. As we stated before, the Hamiltonian flow of $\kplus$ with respect to $\omega_{0}$, is a reparametrization of the geodesic flow on the energy levels and its Hamiltonian vector field is

\begin{displaymath}
X_{\omega_{0}} (\rho) = (\kplus)' (\rho) \, J \, \nabla \rho
\end{displaymath}

where $J$ is an almost complex structure compatible with $\omega_{0}$ such that $J(H(\theta)) = V(\theta)$. Hence, periodic orbits in the level $\rho$ and class $\alpha$ have period

\begin{equation}
\label{eq:period}
T = T (\rho) = \frac{(\kplus)' (\rho)}{\lsub}
\end{equation}

where $\lsub$ is the length of the unique closed geodesic in $\alpha$. 

As the geodesic flow is strongly structurally stable, so is the flow of $\kplus$. This means that there exists an $\bar{\epsilon} > 0$ such that if $\Vert \sigma \Vert _{C^{1}} < \bar{\epsilon}$, then the flows of $\kplus$ with respect to $\omega_{0}$ and $\omega$ are topologically conjugate on the unit energy level, the conjugation being close to the identity. 

Let $\rho = \rho _{0}$ be as in section \ref{constru-k}, i.e. such that $(\kplus)' (\rho_{0}) = - \lsub$. A scaling argument shows that on this energy level, the flows of $\kplus$ with respect to $\omega_{0}$ and $\omega$ are topologically conjugate if $\Vert \sigma \Vert _{C^{1}} < \epsilon$, for $\epsilon = \rho _{0} \bar{\epsilon}$. As $\kplus  \in {\cal{H}} _{V}$, we can make $\epsilon$ independent of $H$ by choosing $\epsilon = \epsilon (V,g) = r_{V} \bar{\epsilon}$, for $r_{V}$ the radius of the fixed neighbourhood $V$. By (\ref{eq:period}), on energy levels near $\rho = \rho_{0}$, periodic orbits of  $X_{\omega_{0}}$ have period $T \approx 1$. As by construction $(\kplus)'' (\rho _{0}) \neq 0$, transversality of (\ref{eq:period}) to $T = 1$ at $\rho = \rho _{0}$ implies that for every magnetic field $\sigma$ with $\Vert \sigma \Vert _{C^{1}} < \epsilon$, there exists $\rho _{\sigma}$ such that the periodic orbits of $X_{\omega}$ in the class $\alpha$ and level $\rho _{\sigma}$ have period $1$. We denote the set of such orbits as ${\cal{P}} _{\omega}(\rho _{\sigma})$. All sets ${\cal{P}} _{\omega_{0}}(\rho)$ are homeomorphic to ${\cal{P}} _{\omega_{0}}(\rho _{0})$ and by Lemma \ref{lema-uno}, this is homeomorphic to $\eseuno$. Then, the structural stability of the flow of $\kplus$ implies that ${\cal{P}} _{\omega}(\rho _{\sigma})$ is homeomorphic to $\eseuno$. By taking a smaller $\epsilon$ if necessary, we can assume that the action of $x \in {\cal{P}} _{\omega}(\rho _{\sigma})$ is still in the interval $[a,b)$ chosen in (\ref{eq:interval}) and that $(\kplus)'' (\rho _{\sigma}) \neq 0$. Then, by Theorem \ref{version}

\begin{equation}
\label{eq:last-step}
\floermag (\kplus; \alpha) \cong H_{\ast} (S^{1}, {\mathbb{Z}}_{2})
\end{equation}

which proves the proposition. $\square$ \\

{\bf Remark} Proposition \ref{iso-mag-st-kplus} still holds when $H$ is not in ${\cal{H}} _{V}$, but with $\epsilon = \epsilon (H,g)$ for $\kminus$. To see this, note that the proof breaks down when trying to find a uniform bound for the magnitude of the magnetic field, as $\epsilon = \rho _{0} \bar{\epsilon} < r_{V} \bar{\epsilon}$, for $\kminus \leq H$, when the support of $H$ is small enough. \\

\begin{figure}[h]

\begin{center}

\psfig{file=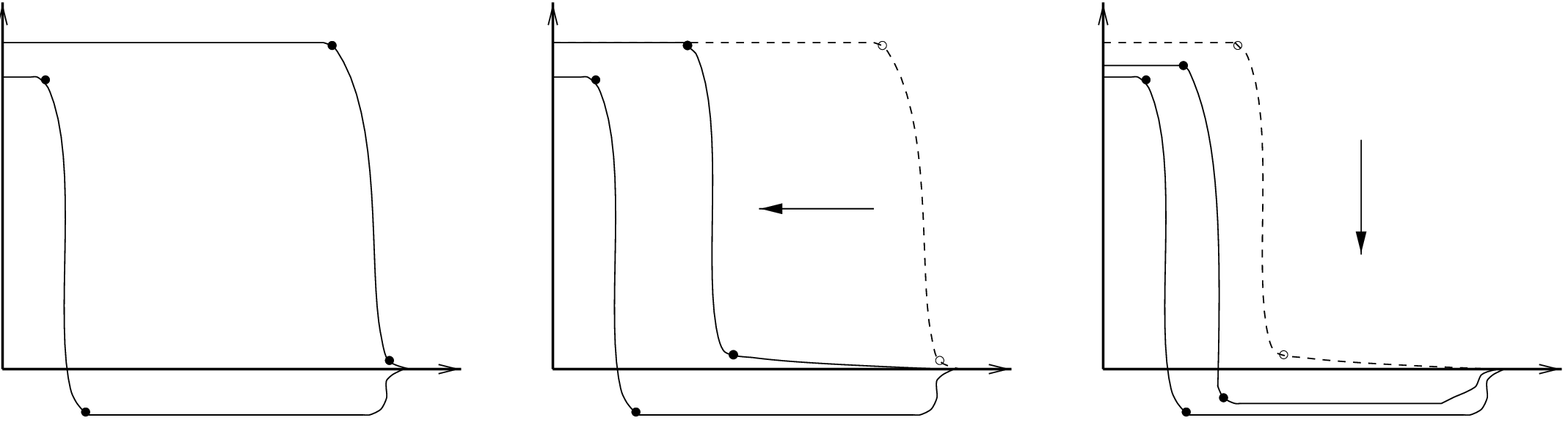, height=12cm, width=4.5cm, angle=270}

\end{center}

\caption{\label{homotopia} Monotone homotopy betwen  $\kplus$ and $\kminus$. }

\end{figure}

{\bf Proof of Theorem \ref{second-theo}} We construct first a monotone homotopy $\{K_{s}\} _{s \in [0,1]}$ connecting $\kplus$ and $\kminus$. Let us take $\kplus$ and, leaving the top of its graph fixed, ``push'' its graph to the left, as shown in Figure \ref{homotopia}. It is clear that the only points at which $(K_{s})' = - \lsub$ are the images by the homotopy of $x^{+}$ and $y^{+}$. Action values for them are decreasing, so $A_{K_{s}} (y^{+})$, where $A_{K_{s}}$ is as in (\ref{eq:functional-small}) goes away from $a$ and $A_{K_{s}} (x^{+})$ decreases to $A_{K_{-}} (x^{-}) > a$.

Then, as seen in Figure \ref{homotopia}, we push the graph obtained in the first step down, completing the homotopy. As before, no new points points at $(K_{s})' = - \lsub$ appear and the action values for $x^{+}$ stay away from $a$. Thus, no new periodic orbits with action equal to $a$ are created. Then, for $\epsilon$ and $\sigma$ as in Proposition \ref{iso-mag-st-kplus}, we conclude that $K_{s} \subset {\cal {H}} ^{[a,b)}$, for $s \in [0,1]$. From Lemmas \ref{lema-floer-1} and \ref{lema-floer-2}, it follows that

\begin{displaymath}
HF _{\omega} ^{[a,b)} (H;\alpha) \cong  H_{\ast} (\eseuno, {\mathbb{Z}}_{2})
\end{displaymath}

which proves the theorem. $\square$

\subsection{Proof of Proposition \ref{propo-1}}
\label{counter}

In this section we prove Proposition \ref{propo-1}. More specifically, we construct a sequence of compactly supported Hamiltonians $K_{n}$, eventually leaving ${\cal{H}} _{V}$ and  with arbitrarily large variation, and a sequence of magnetic terms $\sigma _{n}$ with magnitude tending to zero, such that the periodic orbits of the magnetic flow of $K_{n}$ are all contractible. This shows that Theorem \ref{second-theo} fails if $H$ is not ${\cal{H}} _{V}$. It also shows that the magnitude of the perturbation in Theorem \ref{first-theo} may be arbitrarily small. \\

We describe now the setting in which we work. Let $M$ be a $2$-dimensional Riemannian manifold with negative curvature constantly equal to $-1$ and let $\Omega = dA$ be the area form. We endow the cotangent bundle $\cotgb$ with the symplectic form $\omega _{s} = \omega _{0} + s \pi ^{\ast} \Omega$, where $\omega _{0}$ is the standard symplectic form of $\cotgb$ and $s \in \reals$ is the charge of the magnetic field. We again  identify $\cotgb$ and $TM$ through the metric and we work on $TM$, using the structures introduced in section \ref{neg-curv}.

The twisted geodesic flow is the Hamiltonian flow of $H = \frac{1}{2} \rho ^2$ with respect to $\omega _{s}$. Using the splitting $T_{\theta}TM = H(\theta) \oplus V(\theta)$, we can write the vector field on the energy level $E_{1} = \{ \rho = 1\}$ as

\begin{equation}
\label{eq:horo-vector}
X = J \nabla \rho + s Y _{\Omega} (J \nabla \rho)
\end{equation} 
 
where $J$ is a compatible almost complex structure on $TM$ such that $J(H(\theta)) = V(\theta)$ and $Y_{\Omega}$ is the Lorentz force for $\Omega$. The behaviour of the twisted geodesic flow is qualitatively very different, for different values of $s$. For $s < 1$, it is smoothly diffeomorphic to the geodesic flow. For $s = 1$, it has no closed trajectories. For $s > 1$, all of its orbits are periodic and its projections on $M$ are contractible. \\

Before proceeding with the proof, we give a detailed sketch of it. Let $K = K (\rho)$ be as $\kplus$ in section \ref{constru-k}. By (\ref{eq:slope}), there will be a $1$-periodic orbit in a given non-trivial class $\alpha$ on the energy level $\rho = \rho _{\alpha}$ where $K'(\rho _{\alpha}) = - \lsub$, for $\lsub$ the length of the unique closed geodesic in $\alpha$. As $M$ is closed and negatively curved and $K$ is compactly supported, there are finitely many classes $\alpha _i, i = 1, \cdots, m$ for which the equation

\begin{equation}
\label{eq:ecu-clases}
K' (\rho) = - \lsubi
\end{equation}

holds for some $\rho \in [0,R]$. Note that by condition $d)$ in the construction of $\kplus$ in section \ref{constru-k}, each of the equations (\ref{eq:ecu-clases}) holds at two different energy levels, which we denote by $\rho_{\alpha _i;1} < \rho _{\alpha _i;2}$. Our goal is to construct a magnetic field $\sigma$ such that the Hamiltonian vector field of $H = \frac{1}{2} \rho ^2$ with respect to $\omega = \omega _L + \pi ^{\ast} \sigma$ on each critical level $\rho = \rho _{\alpha _i;j}$ induces dynamics like the one for $s > 1$ in (\ref{eq:horo-vector}), i.e. all orbits on those levels are periodic with contractible projections. Note that this does not contradict Theorem \ref{second-theo}, as the magnitude of $\sigma$ will be large. Then, a scaling argument will give us the desired sequences $K_n$ and $\sigma_n$, with $\Vert \sigma _n \Vert _{C^{1}} \to 0$. \\

{\bf Proof of Proposition \ref{propo-1}} \, Let $K$ be as $\kplus$ in section \ref{constru-k}. We define the magnetic field $\sigma \in \wedge ^2 (M)$ as

\begin{displaymath}
\sigma = c \bar{\rho} \, \Omega
\end{displaymath}

where $\bar{\rho} = \max _{i,j} \{ \rho  _{\alpha _i;j} \}$, $\Omega$ is the area form and $ c > \frac{R}{\bar{\rho}}(1 + \bar{\epsilon})$, where $R$ and $\bar{\epsilon}$ are as in Proposition \ref{iso-mag-st-kplus}. Clearly this makes the $C^1 -$size of $\sigma$ greater than that of the magnetic fields allowed by Theorem \ref{second-theo}. Note that the Lorentz force for $\sigma$ is $Y _{\sigma} = c \bar{\rho} Y_{\Omega}$.

The Hamiltonian vector field of $K$ with respect to $\omega = \omega _L + \pi ^{\ast} \sigma$ at an arbitrary energy level is

\begin{eqnarray*}
\label{eq:vf-in-original-level}
X _{K,\sigma} (\rho) & = & \frac{K'(\rho)}{\rho} (\rho J \nabla \rho + Y_{\sigma} (\rho J \nabla \rho)) \\ & = & \frac{K'(\rho)}{\rho} (\rho J \nabla \rho + c \bar{\rho} Y_{\Omega} (\rho J \nabla \rho)) \\ & = & K'(\rho) (J \nabla \rho + c \bar{\rho} Y_{\Omega} (J \nabla \rho)).
\end{eqnarray*}

Then, at a critical level $\rho = \rho  _{\alpha _i;j}$, we have that 

\begin{equation}
\label{eq:vfincritlevel}
X _{K,\sigma} (\rho  _{\alpha _i;j}) = - \lsubi (J \nabla \rho + c \bar{\rho} Y_{\Omega} (J \nabla \rho)).
\end{equation}

To compare this to (\ref{eq:horo-vector}), we must scale $\rho = \rho  _{\alpha _i;j}$ to $\rho = 1$. As the horizontal part of (\ref{eq:vfincritlevel}) is tangent to the energy level, it is not affected by the scaling, so its pushforward to $\rho = 1$ is 

\begin{displaymath}
X ^{\alpha _{i};j} (1) = - \lsubi (J \nabla \rho + \frac{c \bar{\rho}}{\rho  _{\alpha _i;j}} Y_{\Omega} (J \nabla \rho)).
\end{displaymath}

By our choice of $c$ and $\bar{\rho}$, $\frac{c \bar{\rho}}{\rho  _{\alpha _i;j}} > 1$, so all orbits at all critical levels are periodic with contractible projections.

Now, let $\{ a_{n} \} _{n \in \mathbb{N}}$ be a decreasing sequence, with $a_{n} \to 0$ and $a_{n} \leq 1, \forall n$. We define

\begin{displaymath}
K_{n} (\rho) = K (\rho / a_{n}), \quad \sigma _n = c \bar{\rho _n} \Omega
\end{displaymath}

where $\bar{\rho _n}$ is the analog to $\bar{\rho}$ for $K_n$. Note that for $n$ big enough, $K_{n}$ is not in ${\cal{H}} _{V}$. The argument used for $K$ carries over to $K_n$ and clearly $\Vert \sigma_{n} \Vert _{C^{1}} \to 0$. To complete the proof,  we note that we can make $\min _{M \times \eseuno} K_n$ arbitrarily large by scaling $K_{n}$. $\square$ \\

\end{document}